\documentclass[twoside,openright,11pt,a4paper,epsf]{article}
\usepackage{latexsym}
\usepackage[utf8]{inputenc}
\usepackage{amsmath}
\usepackage{amsthm}
\usepackage{amsfonts}
\usepackage{amssymb}
\usepackage{epsfig}
\usepackage{nicefrac}
\usepackage[all]{xy}
\usepackage{graphicx}

\newcommand{\color}[6]{}
\newcommand{\R}{\mathbb{R}}
\newcommand{\C}{\mathbb{C}}
\newcommand{\D}{\mathbb{D}}
\renewcommand{\P}{\mathbb{P}}

\newcommand{\Z}{\mathbb{Z}}
\newcommand{\Q}{\mathbb{Q}}

\newcommand{\cc}{\mathbb{\mathcal C}}
\newcommand{\cb}{\mathbb{\mathcal B}}
\newcommand{\ca}{\mathcal A}
\newcommand{\ce}{\mathcal E}
\newcommand{\cd}{\mathcal D}
\newcommand{\cv}{\mathcal V}
\newcommand{\cu}{\mathcal U}
\newcommand{\cm}{\mathcal M}
\newcommand{\cn}{\mathcal N}

\newcommand{\cl}{\mathcal{L}}

\newcommand{\nbd}{neighbourhood }
\newcommand{\nbds}{neighbourhoods }
\newcommand{\fonction}[5]
{$$ 
\begin{array}{rcccl}
 #1 & : & #2 & \longrightarrow &#3 \\
    &   & #4 & \longmapsto &#5 
\end{array}
$$}

\newcommand{\priv}{\backslash}
\newcommand{\lra}{\longrightarrow}
\newcommand{\hra}{\hookrightarrow}
\newcommand{\ssi}{\Longleftrightarrow}
\newcommand{\om}{\omega}
\newcommand{\eps}{\varepsilon}
\renewcommand{\phi}{\varphi}
\newcommand{\sdb}{\textnormal{\small SDB}}

\newcommand{\st}{\textnormal{\small st}}

\newcommand{\wdt}[1]{\widetilde{#1}}
\newcommand{\cqfd}{\hfill $\square$ \vspace{0.1cm}\\ }
\newcommand{\sbull}{{\tiny $\bullet$ }}
\newcommand{\ds}{\displaystyle}
\newcommand{\im}{\textnormal{Im}\,}

\renewcommand{\bar}[1]{\overline{#1}}
\newcommand{\pd}{\textnormal{\small PD}}

\newtheorem{definition}{Definition}[section]
\newtheorem{thm}{Theorem}
\newtheorem*{thm*}{Theorem}
\newtheorem{prop}[definition]{Proposition}
\newtheorem{lemma}[definition]{Lemma}

\newtheorem{rk}[definition]{Remark}

\title{ \vspace*{0cm}Singular polarizations and ellipsoid packings.}
\author{Emmanuel Opshtein. \footnote{Partially supported by ANR projects "Floer Power" ANR-08-BLAN-0291-03 and "Symplexe" BLAN06-3-137237}}
\date{}
\begin{document}
\maketitle
\begin{abstract}
We prove in this paper that any $4$-dimensional symplectic manifold is essentially made of finitely many symplectic 
ellipsoids. The key tool is a singular analogue of Donaldson's symplectic hypersurfaces in irrational symplectic 
manifolds.
\end{abstract}
\section{Introduction.}
Donaldson proved in \cite{donaldson} that a symplectic manifold $(M,\om)$ with $\om\in H^2(M,\Z)$ (so-called {\it rational}) 
always admits  a symplectic polarization of large enough degree $k$, that is a symplectic hypersurface Poincaré-dual to 
$k\om$. In \cite{biran3}, Biran showed that these polarizations decompose the manifold into a standard "fat" part and a 
"thin" part which is isotropic in the Kahler case, and which has zero-volume in any case. In \cite{moi3}, it was noticed that 
the standard part of the previous decomposition is itself made of a standard ellipsoid and an object of codimension one. Put together, these results show that rational symplectic manifolds are always covered by {\it one} ellipsoidal Darboux chart up to a negligible set.  This approach is rather satisfactory for $\P^2$ or $(S^2\times S^2,\om\oplus\om)$ where polarizations of low degrees can easily be found. However, as the degree of the polarization becomes larger, the ellipsoid gets more intricate and the codimension-one part more significant. This explosion of degree prevents from getting anything interesting on irrational manifolds. The present work shows however that an analogous result holds in the irrational setting.

\begin{thm}
Any closed $4$-dimensional symplectic manifold has full packing by a finite number of ellipsoids. This number can be 
bounded by  a purely topological quantity :  the dimension of $H^2(M,\R)$.
\end{thm}
This theorem is not really about symplectic embeddings : it does not address the question of how flexible they might be, 
like for instance \cite{schlenk,mcsc,mcduff4,guth,hike}. It rather gives a description of a symplectic manifold as a patchwork of 
euclidean pieces (ellipsoids) whose complexity - if only measured by the number $N$ of pieces - does not really depend on 
the symplectic structure (see also \cite{rusc} for a result in this spirit). Although the bound above is rather loose (for instance when the symplectic form is rational), it can be improved by a closer look at the proof. In fact,
$$
N\leq \min \{\dim V, \; V\subset H^2(M,\Q),\; [\om]\in \textnormal{Span}_\R V\}.
$$
The theorem is a consequence of the following two results. First, Donaldson's construction of polarizations extends to 
irrational symplectic manifolds.
\begin{thm}\label{singpol} For any symplectic manifold $(M,\om)$ there exist symplectic hypersurfaces  $(\Sigma_1,\dots,\Sigma_N)$ with transverse and positive intersections such that 
\begin{equation}\label{eqpol}
[\om]=\sum_{i=1}^N a_i\pd(\Sigma_i), \hspace{,5cm} a_i\in \R^+.
\end{equation}
We can also assume that the classes $\pd(\Sigma_i)$ are independent in $H^2(M,\R)$.
\end{thm}
A family of symplectic hypersurfaces that satisfies (\ref{eqpol}) will be called a 
{\it singular } polarization of $M$. In dimension higher than four, the meaning of "positive intersection" obviously has to be explained, and we refer the reader to section \ref{singpolsec}.  As their classical analogues, singular polarizations can be used to embed ellipsoids, at least in dimension four.
\begin{thm}\label{singdec}
Let $(M^4,\om)$ be a closed symplectic manifold with 
$$
[\om]=\sum_{i=1}^N a_i\pd(\Sigma_i) 
$$
where $\Sigma_i$ are symplectic curves whose pairwise intersections are all positive and whose Poincaré-duals are independant in $H^2(M,\R)$. Then $M$ has a 
full packing by the ellipsoids $\ce(A_i,a_i)$ where $A_i$ denotes the symplectic area of $\Sigma_i$. Precisely,
for all $\eps>0$, there exists an embedding 
$$
\Phi:\amalg \,\ce(A_i-\eps,a_i) \hra M
$$
which admits $(\Sigma_1,\dots,\Sigma_N)$ as supporting surfaces, {\it i.e.} the image of the "horizontal" disc $\{z_2=0\}$
in $\ce(A_i-\eps,a_i)$  covers $\Sigma_i$ up to an area $\eps$.  
\end{thm}
Some remarks are in order. First, a simple computation shows that $M$ is covered by the image of $\Phi$ up to 
a volume of order $\eps$, hence the wording "full packing". Together with theorem \ref{singpol}, it obviously proves the 
basic assertion of this paper. Next, in theorem \ref{singdec}, the curves $\Sigma_i$ are allowed to  have negative self-intersections : the positivity condition only concerns intersections between different curves. As such, it applies for instance naturally in the context of blows-up. It can therefore be used to understand what happens to the ellipsoid 
decomposition in rational sympletic manifolds equipped with polarizations with singularities. It allows in some sense to make 
the desingularization process compatible with Biran's decompositions. Another application concerns symplectic isotopies : 
 the proof of theorem \ref{singdec} goes along the same lines as the proof of Biran's decomposition theorem given in \cite{moi5}, and it extends the range of the method of isotopy developed there. Finally, it may be worth pointing that both the dimension 
 hypothesis and the independance of the Poincaré-duals seem mostly technical, and can be removed at least in some concrete situations ({\it e.g.} $(\P^n,\om_{\textnormal{\tiny FS}})$ with a polarization consisting of two linear hyperplanes is good enough).  

The paper is organized as follows. We first discuss the main idea of the paper through the two easiest examples : the non-singular and the "flat" cases. In section \ref{psdb}, we give a 
local model for a \nbd of a singular polarization, as well as the main properties of this model in terms of Liouville forms. In section \ref{singdecsec}, we prove theorem \ref{singdec}. We then explain the small modifications  to Donaldson's arguments needed to prove the existence of singular polarizations (theorem \ref{singpol}). We finally deal with the applications 
in the last two sections : Biran's decomposition associated to singular curves in section \ref{desingsec} and  isotopies of
 balls in section \ref{isotopysec}.\vspace{-,1cm}
\paragraph{Notations :}We adopt the following  (not so conventional) conventions throughout this paper : \vspace{-,1cm}
\begin{itemize}
\item[-] {\it All} angles will take value in $\R/ \Z$. In other terms an angle $1$ is a full turn in the plane, and the integral of the form $d\theta$ over a circle around the origin in the plane is $1$.\vspace{-,3cm}
\item[-] The standard symplectic form on $\C^n=\R^{2n}$ is $\om_\st:=\sum dr_i^2\wedge d\theta_i$, where $(r_i,\theta_i)$ are polar coordinates on the plane factors. With this convention, the euclidean ball of radius $1$ has capacity $1$. \vspace{-,3cm}
\item[-] A Liouville form $\lambda$ of a symplectic structure $\om$ is a one-form satisfying $\om=\textbf{-} d\lambda$.
The standard Liouville form on the plane is $\lambda_\st:=-r^2d\theta$. \vspace{-,3cm}
\item[-] A symplectic ball or ellipsoid is the image of an euclidean ball or ellipsoid in $\C^n$ by a symplectic embedding. \vspace{-,3cm}
\item[-] The Hopf discs of an euclidean ball in $\C^n$ are its intersections with complex lines. \vspace{-,3cm}
\item[-] $\ce(a,b)$ denotes the $4$-dimensional ellipsoid $\{a^{-1}|z|^2+b^{-1}|w|^2<1\}\subset \C^2(z,w)$. Because of our 
normalizations, its Gromov's width is $\min(a,b)$. 
\end{itemize}

\section{Two easy examples.}
\subsection{The non-singular case.}
In this paragraph we review briefly for self-containedness the result of \cite{moi3} in the setting of smooth polarizations. 
Let $(M,\om)$ be a rational symplectic manifold with a polarization $\Sigma$ of degree $k$. Biran's result states that there is an 
embedding of a symplectic disc bundle $\sdb(\Sigma,k)$ into $M$ which has full volume. This disc bundle can be seen as 
the part of the normal line bundle of $\Sigma$ - denoted by $\cn_\Sigma$  - in $M$ on which the  closed $2$-form $\om_0$ (to be defined soon) is symplectic. The line bundle $\cn_\Sigma$ can be equiped with a hermitian metric and a connection 
form which allow to define a form $\alpha$ on $\cn_\Sigma\priv \cl_0$ satisfying $\alpha_{|F}=d\theta$ and $d\alpha=-k\pi^*\om_{|\Sigma}$. The form $\om_0$ is then simply given in these coordinates by :
$$
\om_0:=\pi^*\om_{|\Sigma}+d(r^2\alpha)=(1-kr^2)\pi^*\om_{|\Sigma}+dr^2\wedge \alpha.
$$
It was proved in \cite{moi3} that the restriction of this disc bundle to a disc of area $A$ in the base is an ellipsoid $\ce(A,1/k)$. 

\begin{lemma}\label{sdb=ell} Let $\pi:\sdb(\Sigma,k)\lra \Sigma$ be the symplectic disc bundle defined above and let 
$D_A$ be a disc of area $A$ in $\Sigma$. Then $(\pi^{-1}(D_A),\om_0)\simeq (\ce(A,1/k),\om_\st)$. 
\end{lemma}
Let us mention an application that was not made completely explicit in \cite{moi3}. It answers a question of McDuff \cite{mcduff4}.
\begin{thm}\label{E(4,1)} There exists a symplectic embedding of $\ce(2,\frac{1}{2})$ into $B^4(1)$.
 \end{thm}
\noindent {\it Proof :} First notice that $\P^2$ has such a full packing because it has a polarization of degree $2$, of area $2$, 
namely a conic. Let us give now an explicit description of a prefered disc bundle over the conic $Q:=\{z_0^2+z_1^2+z_2^2=0\}\subset \P^2$. Since $Q$ is real, it is invariant by conjugacy, and each real projective line intersects 
$Q$ in exactly two distinct conjugated points. Moreover, $\R\P^2$
splits all these lines into two disks of equal area one-half, that contain one of these two points each. The fibers of the disc bundle over the points of $Q$ are precisely these half real lines \cite{biran3}. 
Fix now $z,\bar z\in Q$ and call $d_{z,\bar z}$ the (real) line passing through $z$ and $\bar z$. Consider also a disc $D_Q$ of full 
area which misses $z$ and $\bar z$. The restriction of this symplectic disc bundle to $D_Q$ is an open ellipsoid $\ce(\ca_\om(D_Q),\nicefrac{1}{2})=\ce(2,\nicefrac{1}{2})$. By construction, this ellipsoid does not meet the fibers above $z,\bar z$, so it misses 
the projective line $d_{z,\bar z}$. Since $\P^2\priv d_{z,\bar z}=B^4(1)$, the ellipsoid $\ce(2,\frac{1}{2})$ embeds in fact into $B^4(1)$. \cqfd

Lemma \ref{sdb=ell} serves also to split an ellipsoid into smaller ones. As such, it proved useful to give a natural construction 
of a maximal symplectic packing of $\P^2$ by five balls \cite{moi3}. Let us now mention a far less successful story : looking for such a maximal symplectic packings of $\P^2$ by seven balls (known by \cite{mcpo} to be of radius $r^2=\nicefrac{3}{8}$). 
Using the same idea, one can easily pack $\P^2$ with eight ellipsoids $\ce(\frac{3}{8},\frac{1}{3})$ using a smooth polarization of degree three. These ellipsoids fail to contain the desired balls because $\frac{1}{3}<\frac{3}{8}$. But there are eight of them instead of seven. Notice that one of these ellipsoids can then be split into eight ellipsoids $\ce(\frac{3}{8},\frac{1}{24})$. In this approach, the question would now to be able to glue seven of these eight thin ellipsoids with the seven bigger ones to get seven ellipsoids $\ce(\frac{3}{8},\frac{1}{3}+\frac{1}{24})=\ce(\frac{3}{8},\frac{3}{8})=B^4(\frac{3}{8})$. But this points seems rather hard.
\subsection{The product case.}
Let us discuss now the basic idea of the paper, in the easiest case of a product. Consider the symplectic manifold $M:=(S^2\times S^2,\om\oplus \frac{p}{q}\om)$, where $p,q$ are relatively prime integer. This manifold has a symplectic polarization of degree $q$ which is a smoothening of a curve \fonction{\phi}{S^2}{S^2\times S^2}{z}{(f_1(z),f_2(z)),}
where $f_1,f_2$ are self-maps of $S^2$ of degrees $p$ and $q$ respectively. Over this complicated polarization, 
there is a symplectic ellipsoid $\ce(p+q,\nicefrac{1}{q})$ which cannot be very simple. For instance, when $\nicefrac{p}{q}$ degenerates to an irrational number, Gromov's capacity of the ellipsoid collapses, and nothing remains at the limit. By contrast, there is a much simpler {\it singular} polarization on the homological level given by $(S^2\times \{*\},\{*\}\times S^2)$, which provides a decomposition of $M$ into two ellipsoids $\ce(1,\frac{p}{q})$ in the following way. 

\begin{figure}[h!]
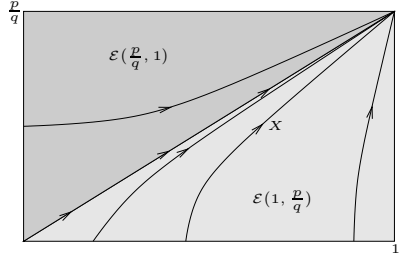

\begin{center} 
\input toricell.pstex_t
\caption{The vector field $X$ on the toric coordinates of $(S^2\times S^2,\om\oplus \frac{p}{q}\om$).}
\label{toricell}
\end{center}
\end{figure}

Put coordinates $\big((r_1,\theta_1),(r_2,\theta_2)\big)$ on $S^2\times S^2$ (remember that $\theta_i\in[0,1[$) with the convention that $\{r_1^2=1\}$ and $\{r_2^2=\nicefrac{p}{q}\}$ is one point ($S^2$ is seen as the one point compactification of the disc of suitable radius). Denote also $\Sigma_1:=S^2\times \{0\}$ and $\Sigma_2:=\{0\}\times S^2$. The symplectic form on $M$ is
$$
\begin{array}{rl}
\om= & dr_1^2\wedge d\theta_1+dr_2^2\wedge d \theta_2 \\
 & =-d\lambda, \hspace{.5cm}\text{ where} \vspace*{,2cm}\\
 \lambda= & (1-r_1^2)d\theta_1+(\frac{p}{q}-r_2^2)d\theta_2.
\end{array}
$$
The Liouville form $\lambda$ is defined on $M\priv (\Sigma_1\cup\Sigma_2)$ and gives rise to a forward complete Liouville vector field, which is easily seen to be 
$$
X=\frac{1-r_1^2}{2r_1}\frac{\partial}{\partial r_1}+\frac{\frac{p}{q}-r_2^2}{2r_2}\frac{\partial}{\partial r_2}.
$$
The action of this vector field is best seen on the toric coordinates $(R_1,R_2):=(r_1^2,r_2^2)$ on $M$, and is shown in figure 
\ref{toricell}.
We  actually see that $X$ is tangent to the line $R_2:=\frac{p}{q}R_1$, so the trajectories of $X$ emanating from 
$\Sigma_1\priv \{(0,0)\}$ and $\Sigma_2\priv \{(0,0)\}$ are respectively $R_1\leq \frac{p}{q}R_2$ and 
$R_2\leq \frac{p}{q}R_1$. These triangles are well-known to be filled by the ellipsoid $\ce(1,p/q)$. Thus we see that 
we get the toric decomposition of $S^2\times S^2$ into two ellipsoids (up to zero volume) out of a data consisting of a singular polarization $(\Sigma_1,\Sigma_2)$ and a Liouville vector field $X$ on the complement of $\Sigma_1\cup\Sigma_2$. This approach provides  much simpler objects (in a geometric sense) than the one giving only one ellipsoid. In particular, both the singular polarization and the embeddings survive the process of degenerating $p/q$ to an irrational. The aim of this note is to understand this simple picture in a  general context.

\section{Plumbed symplectic disc bundles.}\label{psdb}
Let $(M,\om),\Sigma_1,\dots,\Sigma_n$ be as in theorem \ref{singdec}, that is the $\Sigma_i$ are symplectic smooth curves 
with 
$$
[\om]=\sum_1^n a_i\sigma_i,\hspace{1cm} \sigma_i:=\pd(\Sigma_i),\; a_i\geq 0, 
$$
and all intersection points between any two of these curves is positive. Put $\Sigma_i\cap \Sigma_j=:\{(p_{ij}^k)_{k\in[1,l_{ij}]}\}$. With no loss of generality, we can assume that the curves 
are symplectic orthogonal  with respect to $\om$ at each intersection point (such a configuration can be achieved by small local perturbations).
\subsection{Local model near the polarization.}
Decompose first the area form on $\Sigma_i$ as 
$\ds \om_{|\Sigma_i}=\tau_i+\sum_{j,k} \tau_{ij}^k$, where :\vspace{-,2cm}
\begin{itemize}
\item[\sbull] the forms $\tau_{ij}^k$ have supports on small discs $D_{ij}^k$ around $p_{ij}^k$, with total masses $\eps a_j$,\vspace{-,2cm}
\item[\sbull] the form $\tau_i$ has support on the complement $\Sigma_i\priv(\cup D_{ij}^k{}')$ of smaller discs also centered
on $p_{ij}^k$, with total mass
$$
\ca_i^\eps=\ca_\om(\Sigma_i)-\eps\sum_{j\neq i}\Sigma_i\cdot\Sigma_j a_j=a_i\Sigma_i\cdot \Sigma_i+(1-\eps)\sum_{j\neq i}
a_j\Sigma_i\cdot \Sigma_j.\vspace*{-,1cm}
$$
We can also assume that the area of $\tau_i$ on the complement of the discs $D_{ij}^k$ is $\ca_i^{\eps'}$ for $\eps'$ 
slightly smaller than $\eps$. 
\end{itemize}
\begin{figure}[h!]
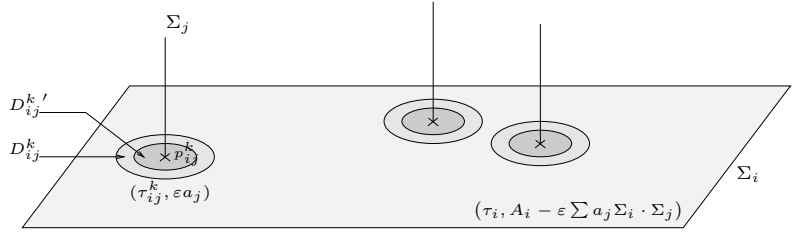

\begin{center} 
\input notations.pstex_t
\caption{Local model near $\Sigma_i$.}
\label{notill}
\end{center}
\end{figure}

Consider now the line bundle $\pi_i:\cl_i\to \Sigma_i$ which is modeled on the (symplectic) normal  bundle of $\Sigma_i$ in 
$M$ - {\it i.e.} they have the same Chern class. Endow this bundle with a hermitian metric, (local) coordinates $(r_i,\theta_i,z)$
and a connection with curvature  $2i\pi\gamma_ia_i^{-1}\pi^*\tau_i$, where 
$$
\frac{1}{\gamma_i}=\frac{a_i^{-1}\ca_{\tau_i}(\Sigma_i)}{\Sigma_i\cdot \Sigma_i}= 1+(1-\eps)\frac{\sum_{j\neq i} a_j\Sigma_j\cdot \Sigma_i}{a_i\Sigma_i\cdot \Sigma_i}
$$
Notice that $\gamma_i$ is negative when $\Sigma_i\cdot \Sigma_i<0$ and vanishes when $\Sigma_i\cdot \Sigma_i=0$. Defining 
the form $\alpha_i$ on $\cl_i$ by asking that its restriction to the fiber is $a_id\theta_i$ and that it vanishes on the 
horizontal planes of the connection, we get a form that checks :
$$
\left \{\begin{array}{l}
\alpha_{i|F}=a_id\theta_i, \\
d\alpha_i=-a_i\gamma_ia_i^{-1}\pi^*\tau_i=-\gamma_i\pi^*\tau_i.
\end{array} \right.
$$
We define now a closed two-form on $\cl_i$ by 
$$
\begin{array}{rl}
\om_i:= & \ds \pi_i^*\tau_i+d(r_i^2\alpha_i)+\sum_{j,k=1}\pi_i^*\tau_{ij}^k, \\
  =& \ds (1-\gamma_ir_i^2)\pi_i^*\tau_i +dr_i^2\wedge \alpha_i+\sum_{k=1}^n \pi_i^*\tau_{ij}^k.
\end{array}
$$
When $\gamma_i$ is non-positive, this form is symplectic on $\cl_i$. But in the positive situation, $\om_i$ is only symplectic on the disc bundle of area $\gamma_i^{-1}$ (on even larger discs over  $D_{ij}^k$). We will denote in the sequel by $\sdb(\cl_i)$ 
the symplectic part of the line bundle. 

A standard Moser argument  shows moreover that there are some \nbds $\cu_i,\cv_i$ of the zero-section $\cl_0$ 
and $\Sigma_i$ respectively which are symplectomorphic. In other terms, there exists an embedding 
$$
\phi_i:(\cu_i,\om_i)\hra(M,\om), \hspace{,5cm} \im \phi_i=\cv_i,\; \phi_i(\cl_0)=\Sigma_i.
$$
For simplicity, we henceforth assume that $\cv_i$ is itself endowed with a fibration (given by $\pi_i\circ \phi_i^{-1}$) and coordinates $(r_i,\theta_i)$. Moreover, since $\Sigma_i$ and 
$\Sigma_j$ are symplectic orthogonal at $p_{ij}^k$, a local symplectomorphism allows to make the fibration structures of $\cv_i$ 
and $\cv_j$ coincide in $\cv_i\cap \cv_j$, namely arranging that $(r_i,\theta_i,r_j,\theta_j)$ provide full coordinate charts 
in $\cv_i\cap\cv_j$, for which the two set of fibres are given by the fibres of $(r_i,\theta_i)$ and $(r_j,\theta_j)$. With such normalization, we can finally assume that 
\begin{equation}\label{tauijk}
\pi_i^*\tau_{ij}^k=a_jdf_{ij}^k(r_j)\wedge d\theta_j, 
\end{equation}
where $f_{ij}^k\equiv \eps$ outside $D_{ij}^k$ and coincides with $r_j^2$ near $p_{ij}^k$.  In some \nbd of this point, we 
therefore have :
$$
\om=a_idr_i^2\wedge d\theta_i+ a_j dr_j^2\wedge d\theta_j.
$$
Let us sum up the above discussion: 
\begin{prop}[Weinstein]\label{weinstein}
There exist \nbds $\cv_i$ of $\Sigma_i$ in $M$ and $\cu_i$ of the zero-section in $\cl_i$ which are identified {\it via} 
diffeomorphisms $\phi_i:\cu_i\to \cv_i$. The expression of the symplectic form in these coordinates is given by 
$$
\phi_i^*\om=\om_i=  \pi_i^*\tau_i+d(r_i^2\alpha_i)+\sum\pi_i^*\tau_{ij}^k,
$$ 
where $\tau_{ij}^k$ has support in $\cv_i\cap \cv_j$ and 
$$
\begin{array}{lll}
\ds\gamma_i=\frac{a_i\Sigma_i\cdot \Sigma_i}{\ca_{\tau_i}(\Sigma_i)}, & 
\left \{\begin{array}{l}
\alpha_{i|F}=a_id\theta_i \\
d\alpha_i=-\gamma_i\pi^*\tau_i
\end{array} \right. , &
\left \{\begin{array}{l}
\pi_i^*\tau_{ij}^k=a_jdf_{ij}^k(r_j)\wedge d\theta_j,\\
f_{ij|{}^cD_{ij}^k}^k\equiv \eps, \hspace{0,5cm} f_{ij}^k=r_j^2 \text{ near } p_{ij}^k
\end{array}. \right. 
\end{array}
$$
Finally, near $p_{ij}^k$, $\om_i=a_idr_i^2\wedge d\theta_i+ a_j dr_j^2\wedge d\theta_j$.
\end{prop}
Otherwise stated,  a \nbd $\cv:=\cup \cv_i$ of the whole polarization is a plumbing of the $\cu_i$ along the 
bidiscs $D_{ij}^k\times D_{ji}^{k'}$ (where $p_{ij}^k=p_{ji}^{k'}$). 

\subsection{Liouville forms on the symplectic disc bundles.}
The symplectic disc bundles $\sdb(\cl_i)$ defined in the previous paragraph come naturally with Liouville forms (recall
they are primitives of the {\it opposite} of the symplectic forms). A more careful analysis - that we perform now - shows that it is 
possible to impose compatibility conditions on these forms, which allow to glue them to get a Liouville form   on $\cv$.
\begin{lemma}\label{liouville}
There is a Liouville form $\lambda_i$ on $\sdb(\cl_i)\priv\big(\cl_0\cup {\pi_i^{-1}(p_{ij}^k)}\big)$ such that 
$\lambda_i=a_i(1-r_i^2)d\theta_i+a_j(1-r_j^2)d\theta_j$  near $p_{ij}^k$. In fact, 
\begin{equation}\label{defliou}
\lambda_i=(1-r_i^2)\alpha_i+(1-\gamma_i)\pi_i^*\lambda_i'+\sum \pi_i^*\lambda_{ij}^k.
\end{equation} 
for well-chosen Liouville forms  $\lambda_i'$, $\lambda_{ij}^k$  for $\tau_i$, $\tau_{ij}^k$ in $\Sigma_i\priv\cup\{p_{ij}^k\}$.
The Liouville form $\lambda_i'$ can however be chosen arbitrarily on any disc compactly supported
 in $\Sigma_i\priv \cup D_{ij}^k$.
\end{lemma}
\noindent {\it Proof :} Consider first any Liouville forms  $\lambda_i',\lambda_{ij}^k$ for $\tau_i,\tau_{ij}^k$ in $\Sigma_i\priv\cup\{p_{ij}^k\}$. Then the one-form defined by (\ref{defliou}) is a Liouville form for $\om_i$. Indeed,
$$
\begin{array}{rl}
d\lambda_i = & -dr_i^2\wedge\alpha_i-(1-r_i^2)\gamma_i\pi_i^*\tau_i -(1-\gamma_i)\pi_i^*\tau_i
-\sum\pi_i^*\tau_{ij}^k \\
 = & -dr_i^2\wedge \alpha_i+(-\gamma_i+\gamma_ir_i^2-1+\gamma_i)\pi_i^*\tau_i-\sum\pi_i^*\tau_{ij}^k \\
  =&-dr_i^2\wedge \alpha_i-(1-\gamma_ir_i^2)\pi_i^*\tau_i-\sum\pi_i^*\tau_{ij}^k \\
   = & -\om_i.
\end{array}
$$
We now need to choose well the forms $\lambda_i'$ and $\lambda_{ij}^k$. Define first $\lambda_{ij}^k$ by
$$
\lambda_{ij}^k:= a_j(\eps- f_{ij}^k(r_j))d\theta_j,
$$ 
and recall that by definition of $f_{ij}^k$, it vanishes identically outside $D_{ij}^k$. In order to define $\lambda_i'$,
notice that $\tau_i$ has support in $\Sigma_i\priv {D_{ij}^k}'$ and $(1-\gamma_i)\ca_{\tau_i}(\Sigma_i)=(1-\eps)\sum_{j\neq i} a_j\Sigma_i\cdot \Sigma_j$. Therefore, there exists a Liouville form $\lambda_i'$ of $\tau_i$ such that 
$$
(1-\gamma_i)\lambda_i'=(1-\eps) a_jd\theta_j \hspace{1cm}\text{ near } p_{ij}^k.
$$
It is moreover obvious  that this condition is compatible with any requirement 
on $\lambda_i'$ on a disc compactly supported in $\Sigma_i\priv \cup D_{ij}^k$.
Putting all this together, we get the following expression for $\lambda_i$ in the \nbd of $p_{ij}^k$ :
$$
\hspace*{2cm}\begin{array}{rlr}
\lambda_i=  & (1-r_i^2)\alpha_i+(1-\eps) a_jd\theta_j+ a_j(\eps-r_j^2)d\theta_j &\\
 = &a_i(1-r_i^2)d\theta_i+a_j(1-r_j^2)d\theta_j. & \hspace*{1,8cm} \square
\end{array}
$$
Recall that a Liouville form $\lambda$ gives rise to a vector field $X_\lambda$ - called Liouville - by symplectic 
duality : $\iota_{X_\lambda}\om=\lambda$. This vector field has the property of contracting the symplectic form : 
$\Phi_{X_\lambda}^{t*}\om=e^{-t}\om$. Thanks to the cautious choices we made until now,  both the sets 
of Liouville forms $(\lambda_i)$ and  vector fields $(X_{\lambda_i})$ glue together to well-defined objects on $\cv\priv (\cup \Sigma_i)$. 
\begin{lemma}\label{gluliou}
The formulas 
$$
\left\{\begin{array}{ll}
\lambda_{|\cv_i}&:=\phi_{i*}\lambda_i \\
X_{\lambda |\cv_i} & := \phi_{i*}X_{\lambda_i}
\end{array} \right.
$$
define a Liouville form and its associated Liouville vector fields on $\cv\priv \cup \Sigma_i$. Moreover, the vector field 
$X_\lambda$ points outside $\cv$ if this \nbd is well-chosen.
\end{lemma}
\noindent {\it Proof :} The first point is an obvious consequence of the previous lemma because $\lambda_i=\lambda_j$ 
near $p_{ij}^k$. The second statement is a straightforward consequence from the fact that each $X_{\lambda_i}$ points
outside the zero-section on $\cl_i$, and this is a simple computation :
$$
\begin{array}{rcl}
\ds\om_i\left(X_{\lambda_i},\frac{\partial}{\partial \theta_i}\right)= & \ds dr_i^2\wedge \alpha_i\left(X_{\lambda_i},\frac{\partial}{\partial \theta_i}\right) & =a_idr_i^2\big(X_{\lambda_i}\big) \\
 = &\ds \lambda_i\left(\frac{\partial}{\partial \theta_i}\right)&=(1-r_i^2)a_i \text{ (see (\ref{defliou})).}
\end{array}
$$
Thus $dr_i^2(X_{\lambda_i})=1-r_i^2>0$ near the zero-section $\{r_i=0\}$.\cqfd

The following lemma gives a nice expression of the Liouville vector fields 
associated to the forms $\lambda_i$ defined above. In the statement, the disc $\cd_A$  should be thought of as a disc of 
$\Sigma_i\priv \cup D_{ij}^k$ of approximately full area. 
\begin{lemma}\label{cse} Consider the trivial disc bundle $\pi:\cd_A\times \D_{\gamma^{-1}}\lra\cd_A$ (or $\cd_A\times \C$ if $\gamma<0$) over a disc in $\C$, with polar coordinates $(r,\theta)$ and $(\rho,\zeta)$ on $\D_{\gamma^{-1}}$ and $\cd_A$ 
respectively. Equip this bundle with the symplectic structure  $\om:=\pi^*\om_{\textnormal{st}}+d(r^2\alpha)$, where 
$\alpha_{|\{x\}\times \D}=a d\theta$ and $d\alpha=-\gamma\pi^*\om_{\textnormal{st}}$. Let $\lambda$ be a Liouville form for 
$\om$ defined by 
$$
\lambda=(1-r^2)\alpha+(1-\gamma)\pi^*\lambda_{\textnormal{st}}.
$$
and $X_\lambda$ its associated vector field. Then 
\begin{itemize}
\item[i)] 
$$
X_\lambda=\frac{1-r^2}{2r}\frac{\partial}{\partial r}-\frac{1-\gamma}{1-\gamma r^2}\frac{\rho}{2}\frac{\partial}{\partial \rho};
$$
\item[ii)] there exists a smooth function $h:\cd_A\lra
  \R$ such that the map \fonction{\Phi}{(\cd_A\times 
  \D_{\gamma^{-1}},\om)}{\big(\ce(A,a\gamma^{-1}),\om_{\textnormal{st}}\big)}{(z,w)}{(z',w')=(\sqrt{1-\gamma|w|^2}z,\sqrt a
   e^{i\,h(z)}w)} is a symplectomorphism (when $\gamma$ is negative, $\ce(A,a\gamma^{-1})$ is an hyperboloid 
   rather than an ellipsoid);
\item[iii)] setting $R':={r'}^2=|w'|^2$ and ${\cal P}':=\rho'^2=|z'|^2$, 
$$
\Phi_*X_\lambda=(a-R')\frac{\partial}{\partial R'}-{\cal P}'\frac{\partial}{\partial {\cal P}'}+*\frac{\partial}{\partial \theta'}.
$$
\end{itemize}
\end{lemma}
\noindent {\it Proof :} The point ii) is word for word the same statement and same proof than lemma 2.1 in \cite{moi3}.
It is an easy computation, which we do not repeat here. The point i) is a simple verification. Write $\om=(1-\gamma r^2)d\rho^2\wedge d\zeta+dr^2\wedge \alpha$ and compute :
$$
\begin{array}{rl}
\om(\frac{1-r^2}{2r}\frac{\partial}{\partial r}-\frac{1-\gamma}{1-\gamma r^2}\frac{\rho}{2}\frac{\partial}{\partial \rho},\cdot) = &
(1-r^2)dr\wedge \alpha(\frac{\partial}{\partial r},\cdot)-(1-\gamma)\rho^2d\rho\wedge d\zeta(\frac{\partial}{\partial \rho},\cdot)\\
 = & (1-r^2)\alpha-(1-\gamma)\rho^2d\zeta\\
  = & \lambda.
\end{array}
$$
For iii), express first $\Phi$ in the good coordinates $\Phi({\cal P},\zeta,R,\theta)=({\cal P}',\zeta',R',\theta')$ :
\begin{equation}\label{eq1}
{\cal P}'=(1-\gamma R){\cal P},\hspace{,5cm} R'=aR, \hspace{,5cm} \zeta'=\zeta \hspace{,5cm} \theta'=\theta+h({\cal P},\zeta).
\end{equation}
Then,
\begin{equation}\label{eq2}
\left\{\begin{array}{l}
\ds \Phi_*\frac{\partial}{\partial R}=-\gamma {\cal P} \frac{\partial}{\partial {\cal P}'}+a\frac{\partial}{\partial R'} \\
\ds \Phi_*\frac{\partial}{\partial {\cal P}} =(1-\gamma R)\frac{\partial}{\partial {\cal P}'}+*\frac{\partial}{\partial \theta'}.
\end{array} \right.
\end{equation}
Taking (\ref{eq1}) and (\ref{eq2}) into account, we therefore get  :
$$
\hspace*{,3cm}\begin{array}{rlr}
\Phi_*X_\lambda= & \ds\Phi_*\big((1-R)\frac{\partial}{\partial R}-\frac{1-\gamma}{1-\gamma R}{\cal P}\frac{\partial}{\partial {\cal P}}\big)& \\
\overset{\tiny(\ref{eq2})}{=}& \ds(1-R)\big[-\gamma {\cal P} \frac{\partial}{\partial {\cal P}'}+a\frac{\partial}{\partial R'}\big]-(1-\gamma) {\cal P}\frac{\partial}{\partial  {\cal P}'}+*\frac{\partial}{\partial\theta'} & \\
 = &  \ds a(1-R)\frac{\partial}{\partial R'}-\big[-\gamma {\cal P}+\gamma {\cal P} R- {\cal P}+\gamma {\cal P}\big]\frac{\partial}{\partial  {\cal P}'}+*\frac{\partial}{\partial \theta'} & \\
\overset{\tiny(\ref{eq1})}{=}&\ds (a-R')\frac{\partial}{\partial R'}- {\cal P}'\frac{\partial}{\partial  {\cal P}'}+*\frac{\partial}{\partial \theta'}. &\hspace*{,8cm} \square
\end{array}
$$

\subsection{Ellipsoids in the standard bundles.}\label{embell}
The ellipsoids of theorem \ref{singdec} naturally arise from $\sdb(\cl_i)$ as the set of points that can be reached by flowing
out of a disc in $\Sigma_i\priv\cup D_{ij}^k$ along the Liouville vector field. Precisely,
\begin{prop}\label{ellinsdb}
Let $\cd_{A_i-\delta}$ be a disc of symplectic area $A_i-\delta$ in $\Sigma_i\priv \cup D_{ij}^k$, viewed as the zero-section of
$\sdb(\cl_i)$. Then, if the form $\lambda_i$ is well-choosen on $\cd_{A_i-\delta}$, the basin of attraction of this disc, defined as 
$$
\cb_i:=\left\{p\in \sdb(\cl_i)\;|\; \exists t\in \R^+, \; \Phi^{-t}_{X_{\lambda_i}}(p)\in \cd_{A_i-\delta}\right\}
$$
is symplectomorphic to the ellipsoid $\ce(A_i-\delta,a_i)$.
\end{prop}
\noindent {\it Proof :} Since $\cd_{A_i-\delta}$ is contained in $\Sigma_i\priv \cup D_{ij}^k$, the symplectic form on 
the restriction of $\sdb(\cl_i)$ to $\cd_{A_i-\delta}$ is exactly of the form of lemma \ref{cse} :
$$
\left\{\begin{array}{l}
\om_i=\pi_i^*\tau_i+d(r_i^2\alpha_i) \\
\lambda_i=(1-r^2)\alpha_i+(1-\gamma_i)\pi_i^*\lambda_i',\hspace{,3cm} d\lambda_i'=-\tau_i.
\end{array}\right.
$$
Provided $\lambda_i'$ corresponds also to the Liouville form called "standard" in this lemma (which can always be achieved because $\lambda'_i$ can be any Liouville form on $\cd_A$ by lemma \ref{liouville}), 
it provides a symplectic embedding $\Phi:\big(\pi_i^{-1}(\cd_{A_i-\delta}),\om_i\big)\hra \big(\C^2,\om_\textnormal{st}\big)$. This map sends the set $\cb_i$ to 
$$
\Phi(\cb_i)=\left\{p\in \C^2\; |\; \exists t\in \R^+,\; \Phi^{-t}_{\Phi_*X_{\lambda_i}}(p)\in \D_{A_i-\delta}\times \{0\}\right\}.
$$
By lemma \ref{cse} iii), if $(z,w)$ are coordinates on $\C^2$, $ {\cal P}=|z|^2$ and $R=|w|^2$, the differential equation 
associated to $\Phi_*X_{\lambda_i}$ is 
$$
\left\{\begin{array}{l} \dot R=a_i-R \\ \dot {\cal P}=-{\cal P} \end{array}\right. , \text{ with solutions }
\left\{\begin{array}{l} R(p,t)=a_i-c_1(p)e^{-t} \\ {\cal P}(p,t)=-c_2(p)e^{-t} \end{array}\right. .
$$
Now $\Phi(\cb_i)$ is the set of points $p\in \C^2$ that verify :
\begin{equation}\tag{$*$}
R(p,t_0)=0\Longrightarrow {\cal P}(p,t_0)<A_i-\delta.
\end{equation}
An easy computation shows that ${\cal P}(p,t_0)=\frac{c_2(p)}{c_1(p)}a_i$, so that ($*$) writes 
$c_2(p)a_i\leq c_1(p)(A_i-\delta)$. This in turn means 
$$
\begin{array}{rll}
 \ds R(p)=a_i-\frac{c_1(p)}{c_2(p)}{\cal P}(p)\leq a_i-\frac{a_i}{A_i-\delta}{\cal P}(p) \ssi & \ds \frac{R(p)}{a_i}+\frac{{\cal P}(p)}{A_i-\delta}\leq 1 & \\
 \ssi & p\in \ce(A_i-\delta,a_i). & \square
\end{array}
$$

We conclude this paragraph by noting that this ellipsoid is contained in the part of the bundle above the disc $\cd_{A_i-\delta}$
simply because of the formula i) of lemma \ref{cse}. Indeed, since $\gamma<1$, the "horizontal" part 
$$
-\frac{1-\gamma}{1-\gamma r^2}\frac{\rho}{2}\frac{\partial}{\partial \rho}
$$
of the vector field $X_\lambda$ above $\cd_{A_i-\delta}$ points inside $\cd_{A_i-\delta}$. 
\begin{rk}\label{localell} The set $\cb_i$ lies inside $\pi_i^{-1}(\cd_{A_i-\delta})$.
\end{rk}
 
\subsection{Variations of the Liouville forms.}\label{variations}
Liouville forms are never unique : they can always be modified by adding a closed one-form. In the previous paragraphs, we needed to impose several compatibility conditions for the Liouville forms, namely fix them on discs $\cd_{A_i-\delta}, (D_{ij}^k)$.
These requirements only rigidify slightly the situation but still leaves a lot of freedom, which will be fully needed in the proof of theorem \ref{singdec}. Precisely, we will need the following set of objects :
\begin{itemize}
\item[\sbull] A family $\vartheta:=(\vartheta_i)$ of closed one-forms  on $\Sigma_i$ which {\it vanish identically} on all 
the $D_{ij}^k$ and $\cd_{A_i-\delta}$. Notice that all homological classes in $H^1_\textnormal{DR}(\Sigma_i)$ have such representatives. 
\item[\sbull] A family $\lambda_\vartheta:=(\lambda_i+\pi_i^*\vartheta_i)$ of Liouville forms on $\sdb(\cl_i)$. 
\end{itemize} 
These forms obviously  satisfy the same compatibility conditions as the $(\lambda_i)$, {\it i.e.} they give rise to a well-defined Liouville form still denoted $\lambda_\vartheta$ on $\cv\priv \cup\Sigma_i$. Moreover, since $\lambda_\vartheta=\lambda$
in $\cd_{A_i-\delta}$ (and therefore in $\pi_i^{-1}(\cd_{A_i-\delta})$), the remark  \ref{localell} ensures that proposition \ref{ellinsdb} holds when $\lambda$ is replaced by $\lambda_\vartheta$. Finally, since $\lambda_\vartheta$ differs from 
$\lambda$ only by a pull-back by $\pi_i$, the radial component of its Liouville vector field does not change : it still moves away from the zero-section, so that lemma \ref{gluliou} also holds for $\lambda_\vartheta$. 

\section{Proof of theorem \ref{singdec}.}\label{singdecsec}
We adopt in this paragraph all conventions, notations and results of section \ref{psdb} . The core lemma is now the following :
\begin{lemma}\label{core}
There exists a family of one-forms $(\vartheta_i)$ on $\Sigma_i$ which vanish identically on $\cd_{A_i-\delta}$ and $D_{ij}^k$
such that the form $\lambda_\vartheta$ defined on $\cv\priv \cup \Sigma_i$ extends to a Liouville form $\beta$ on $M\priv \cup \Sigma_i$.
 \end{lemma} 
Let us first explain quickly why theorem \ref{singdec} is a direct consequence of this lemma. Since $M$ is compact and 
$X_\vartheta:=X_{\lambda_\vartheta}$ points outside the $\Sigma_i$, it defines a forward-complete vector field on $M\priv \cup \Sigma_i$. Therefore, and since $\beta$ is really an extension of $\lambda_\vartheta$, the elementary dynamical procedure 
that consists in extending the local symplectic embeddings $\phi_i:\cu_i\hra \cv_i$ (given by prop \ref{weinstein}) by
\fonction{\Phi_i}{\sdb(\cl_i)}{M}{x}{\left \{\begin{array}{cl} \phi_i(x) & \text{ if } x\in \cu_i \\ \Phi_{X_\beta}^\tau\circ \phi_i\circ\Phi_{X_{\vartheta_i}}^{-\tau}(x) & \text{ if } \Phi_{X_{\vartheta_i}}^{-\tau}(x)\in \cu_i \end{array}\right.}
provides symplectic embeddings $\Phi_i$ which overlap, but clearly not on the sets $\cb_i$. We therefore have an embedding 
$\Phi:\cup \cb_i\hra M$ which is the desired ellipsoid packing by proposition \ref{ellinsdb}.\cqfd

\noindent {\it Proof of lemma \ref{core} :} First observe that by definition of the curves $\Sigma_i$, the symplectic form $\om$
vanishes on any cycle of $M\priv \cup\Sigma_i$, so it is exact on $M\priv \cup\Sigma_i$, and we can pick a Liouville form 
$\beta$ for $\om$ on this set. In $\cv\priv \cup\Sigma_i$, the difference $\beta-\lambda$ is closed. If it is moreover exact, 
the lemma follows because any extension of the function $h$ defined by $\beta-\lambda=dh$ gives an extension $\beta-dh$
of $\lambda$ to the whole of $M\priv \cup \Sigma_i$. We explain now that although this difference may well not be exact, we can find a "correction" closed one-form  $ \vartheta$ as in paragraph \ref{variations} such that $\beta-\lambda_\vartheta$ is exact. 
To understand this point, consider a family $\{\gamma_\eps^i,\gamma_l^i\}_{i,k}$ for the one-dimensional homology 
of $\cv\priv \cup \Sigma_i$, where $\gamma_\eps^i$ is the small loop around $\Sigma_i$ (contained in a fibre of $\pi_i$ 
and defined by the equation $r_i=\eps$) and the $\gamma_l^i$ are $\pi_i$-lifts of simple closed loops ${\gamma_l^i}'$ in 
$\Sigma_i$ which span $H_1(\Sigma_i)$.

We first prove that $f_i:=[\beta-\lambda](\gamma_\eps^i)$ vanishes for all $i$. Since $\lambda(\gamma_\eps^i)$ tends 
to $a_i$ when $\eps$ goes to zero, $\beta(\gamma_\eps^i)$ also has a limit, $a_i+f_i$. Consider now 
a two-cycle $C\in H_2(M)$ and perturb it so that it becomes transverse to the curves $\Sigma_i$. Then since $d\beta=-\om$,
we have :
$$
\begin{array}{rl}
\ds \int_C\om=& \ds\sum_i \lim \beta(\gamma_\eps^i) C\cdot\Sigma_i \\
= & \ds \sum(a_i+f_i) C\cdot \Sigma_i \\
= & \ds \om([C]) + \sum f_i\Sigma_i\cdot C.
\end{array}
$$
Thus, $\sum f_i\pd(\Sigma_i)$ vanishes in $H^2(M,\R)$ which implies the vanishing of each $f_i$ by the independence hypothesis. Notice that we only use this hypothesis at this point, so when the $\gamma_\eps^i$ are contractible for instance, 
the independance is not needed.

Define now $\vartheta$ by requiring that $\vartheta_i\cdot [{\gamma_l^i}']=\int_{\gamma_l^i} \beta-\lambda$. Provided 
that we were cautious to take ${\gamma_l^i}'$ with no intersection with $\cd_{A_i-\delta}$ and $D_{ij}^k$, we can even require 
$\vartheta_i$ to vanish on these discs. Then a simple computation (explicitly made in \cite{moi5}) shows that $\beta-\lambda_\vartheta$ vanishes on each class $[\gamma_l^i]\in H_1(\cv\priv \cup\Sigma_i)$. Moreover, since $\lambda_\vartheta=\lambda+\pi^*\vartheta$, its values on the loops $\gamma_\eps^i$ remain unchanged, so that $[\lambda_\vartheta-\beta](\gamma_\eps^i)=0$ also. The form $\beta-\lambda_\vartheta$ has therefore no period in $\cv\priv \Sigma_i$, so it is exact. \hfill $\square$

\section{Existence of singular polarizations.}\label{singpolsec}
We now prove theorem \ref{singpol}, which asserts that singular polarizations always exist. Let us fix a symplectic manifold $(M,\om)$. We have to find a decomposition of the cohomolgy class of the symplectic form into a sum of 
Poincaré-dual of symplectic hypersurfaces $\Sigma_i$ which intersect transversally and positively. 
Of course, positive intersections is well-defined only in dimension four. In higher dimensions, we model the definition to complex manifolds.
\begin{definition} 
Symplectic submanifolds $\Sigma_1,\dots,\Sigma_k$  of $(M^{2n},\om)$ are said to 
intersect transversely and positively if   all intersections $\Sigma_{j_1}\cap\dots\cap \Sigma_{j_p}$ between $p$ of these submanifolds are transverse, generic  and symplectic.
\end{definition}

\subsection{Proof of theorem \ref{singpol}.}
First notice that the assumption on the independence of the classes $\sigma_i:=\pd(\Sigma_i)$ can be freely removed. 
Indeed, if there is a decomposition of $[\om]$ as in theorem \ref{singpol} with a linear relation 
$\sum \lambda_i\sigma_i=0$, $\lambda_i\in \R$, assume - after maybe changing the indices - that $|a_N\lambda_N^{-1}|\leq |a_i\lambda_i^{-1}|$ for all $i$ (in particular $\lambda_N\neq 0$). Then, writing 
$\sigma_N=-\sum_{i\leq N-1} \frac{\lambda_i}{\lambda_N}\sigma_i$, we get 
$$
[\om]=\sum_{i\leq N-1} (a_i-\frac{\lambda_i}{\lambda_{N}}a_N)\sigma_i=\sum_{i\leq N-1} a_i'\sigma_i \text{ and } a_i'\geq 0.
$$

Let $(M,\om)$ be our symplectic manifold and write $\om=\sum_1^k b_i\sigma_i$, where $[\sigma_i]\in H^2(M,\Z)$. The real vector $b=(b_i)$ is a barycenter of  $N$ rational vectors nearby (at most $\dim H^2(M,\R)+1$), that is for any small $\eps>0$ we have 
$$
b=\sum_1^{N}\lambda_jb^j, \hspace{,3cm} \sum_1^{N}\lambda_j=1, \hspace{,3cm} \Vert b-b^j\Vert<\eps,\hspace{,3cm} b^j\in \Q^N.
$$
Thus, 
$$
\begin{array}{rll}
\om=  & \ds \sum_{i=1}^k\big(\sum_{j=1}^{N}\lambda_jb_i^j\big)\sigma_i &\\
 =& \ds \sum_{j=1}^{N} \lambda_j\big(\sum_{i=1}^k b_i^j\sigma_i\big) & \ds  =\sum_{j=1}^{N} \lambda_j\om_j,
\end{array}
$$
where $|\om-\om_j|<\eps$ and $\om_j\in H^2(M,\Q)$. If $\eps$ is small enough, the forms $\om_j$ are symplectic, so by a result of Donaldson, there are $\om_j$-symplectic hypersurfaces $(\Sigma_1,\dots,\Sigma_{N})$ and positive integers $k_1,\dots,k_{N}$ such that $\pd(\Sigma_j)=k_j\om_j$. Thus,
$$
[\om]=\sum_{j=1}^{N} a_j\pd(\Sigma_j), \hspace{,5cm} a_j=\frac{\lambda_j}{k_j}\in \R^+.
$$
Recall at this point that the $\Sigma_j$ are known to be $\om_j$-symplectic because they are almost $J_j$-holomorphic 
for an $\om_j$-compatible almost-complex  structures. Now if the $(\om_j)$ are close enough to $\om$, all the $J_j$  tame $\om$, 
so the $\Sigma_j$ are also $\om$-symplectic.
What remains to show is that the $\Sigma_j$ can be required to meet positively and transversely. This point proceeds
from the following theorem, which we only state in dimension $4$ for simplicity, but whose generalization to higher dimensions is straightforward.
\begin{thm}\label{da} Let $(M^4,\om,J,g)$ be a symplectic manifold with a compatible almost-complex structure and its associated 
metric. Let $\om_i$ be rational symplectic forms on $M$ close to $\om$, with $(\om_i,g)$-compatible almost complex structures 
$J_i$.  Let $\cl_i\to M$ be a hermitian line bundle endowed with a  connection of curvature $2i\pi q\om_i$ ($q$ being such that 
$q\om_i\in H^2(M,\Z)$ for all $i$). Denoting $g_k:=kng$, there exist sequences of sections $s_i=(s_i^k)$ of 
$\cl_i^{\otimes k}$ such that :
\begin{itemize}
\item[i)] $s_i$ is approximately $J_i$-holomorphic, {\it i.e.} :
$$
|s_i^k|_{g_k,\cc^1}\leq C,\hspace{,3cm} |\bar \partial_{J_i} s_i^k|_{\cc^1,g_k}\leq C/\sqrt k \hspace{,3cm} \text{ for large } k,
$$
\item[ii)] $s_i$ is $\eta$-transverse to $0$, {\it i.e.} $|s_i^k|\leq \eta\Rightarrow |\partial_{J_i} s_i^k|\geq \eta$,
\item[iii)] for all $(i,j)$, the sequence of sections $(s_i,s_j)$ of $\cl_i^{\otimes k}\oplus \cl_j^{\otimes k}\to M$ is $\eta$-transverse 
to $0$, {\it i.e.} :
$$
\begin{array}{rl}
\forall p\in M, |(s_i^k,s_j^k)|<\eta \Longrightarrow & (\partial_{ J_i} s_i^k,\partial_{J_j} s_j^k)\in \cl(T_pM,\C^2) \text{ has a right  } \\
 &  \text{ inverse of $g_k$-norm less than } \eta^{-1}, 
\end{array}
$$ 
\item[iv)] For all $(i,j,l)$, the section $(s_i,s_j,s_l)$ of $\cl_i^{\otimes k}\oplus \cl_j^{\otimes k}\oplus \cl_l^{\otimes k}\to M$
is $\eta$-transverse, {\it i.e.} it has norm at least $\eta$.
\end{itemize} 
\end{thm}

It is important to notice that  in the theorem above, everything concerns {\it sequences } of sections, the norm involving $s_i^k$ 
is {\it always} $g_k:=kng$, and the constants $C$ and $\eta$ depend neither on $k$ {\it nor on the choice of the symplectic structures } $\om_i$ provided they are on a small \nbd of $\om$. This theorem is unfortunately not a 
formal consequence of theorems already stated  in Donaldson and Auroux's papers on the subject because  the bundles 
we consider are of the form $\cl_i^{\otimes k}\oplus \cl_j^{\otimes k}$ instead of $E\otimes \cl^{\otimes k}$. Although 
it however follows from their proofs themselves, we choose in this paper to review (once again) Donaldson's technique (with Auroux's contributions) and to include in the discussion the small modifications we need to make in our setting. This will be done in the next paragraph. We now explain why theorem \ref{da} indeed implies theorem  \ref{singpol}.

As we already noticed,  the vanishing sets of $s_i^{k_0}$ for $k_0\gg 1$ (which we denote $s_i$ in the sequel since $k_0$ is fixed) give $\om$-symplectic hypersurfaces $\Sigma_i\subset M$ such that 
$$
[\om]=\sum a_i\pd(\Sigma_i).
$$
We need to understand that the transversality conditions iii) and iv) implies transversality and positivity of the intersections  between $\Sigma_i$ and $\Sigma_j$. First, condition iv) obviously implies that the intersections are simple : they never involve more than two branches. Let now $p\in \Sigma_i\cap \Sigma_j$, that is $s_i(p)=s_j(p)=0$.  In order to show that 
the intersections between $\Sigma_i$ and $\Sigma_j$ is positive at $p$, we make the following two observations :
\begin{enumerate}
\item $T_p\Sigma_i$ and $T_p\Sigma_j$ are very close (for $k_0$ large enough) to $J_{i/j}$-holomorphic hyperplanes (=lines in the 4-dimensional situation) $\Pi_i,\Pi_j$ 
in $T_pM$.
\item The angle between $\Pi_i$ and $\Pi_j$ is bounded from below by some constant  $C(\eta)$ depending neither on $k$ nor on the symplectic structures $\om_i$. 
\end{enumerate} 
Taking the complex structures $J_i$ close to $J$ by an amount $\eps\ll C(\eta)$, we therefore find that $T_p\Sigma_i$ and 
$T_p\Sigma_j$ are $\eps$-close to $J$-holomorphic lines which form an angle approximately $C(\eta)$-large. Since two $J$-holomorphic lines intersect positively when they are disjoint, we conclude that $T_p\Sigma_i\cap T_p\Sigma_j$ is a positive transverse intersection. Point (1) is very classical and is at the core of Donaldson's proof. Point (2) is only slightly more involved linear algebra done in \cite{auroux1}. Let us prove them anyway.  

Write $ds_i(p)=u_i+\eps_i$ where $u_i=\partial_{J_i} s_i(p)$, $\eps_i=\bar \partial_{J_i} s_i(p)$. Then by i) $|\eps_i|\ll 1$ if $k_0$ is large enough (recall that $|\cdot|$ means $|\cdot|_{g_{k_0}}$), while $|u_i|\geq \eta$ by ii) and 
$$
(u_i,u_j):T_pM\lra \C^2
$$ 
is invertible (recall that $\dim_\R T_pM=4=\dim\C^2$ so right-invertible means invertible) with inverse $R$ of norm less 
than $\eta^{-1}$ by iii). 

To understand ($1$), notice that $T_p\Sigma_i=\ker ds_i(p)=\ker (u_i+\eps_i)$ and  consider  a unitary vector 
$x\in T_p\Sigma_i$  decomposed as $x_0+\tau$ with $x_0\in \ker u_i$ and $\tau \perp_{g_{k_0}} \ker u_i$. Then, 
$$
u_i+\eps_i(x)=0=u_i(\tau)+\eps_i(x),
$$ 
so $u_i(\tau)=-\eps_i(x)$. Taking into account that $\tau\in (\ker u_i)^\perp$ we know that $|u_i(\tau)|=|u_i||\tau|$, so 
$$
|\tau|\leq \frac{|\eps_i|}{|u_i|}\ll 1.
$$
Therefore, $x$ is close to a unitary vector in $\Pi_i:=\ker u_i=\ker \partial s_i(p)$, so $T_p\Sigma_i$ is close (in the angle sense)
to the $J_i$-holomorphic hyperplane $\Pi_i$. In order to estimate the angles between $\Pi_i$ and $\Pi_j$, put
$$
\begin{array}{rl}
\kappa:= & \min\{|\langle x,y\rangle|, x\in \ker u_i,y\in \ker u_j, |x|=|y|=1\} \\
 = & \min\{|\pi_i(y)|, y\in \ker u_j, |y|=1\},
\end{array}
$$ 
where $\pi_i$ stands for the $g_{k_0}$-orthogonal projection on $\ker u_i$. Then, $\kappa=\cos \theta$ where 
$\theta$ is the angle between $\Pi_i$ and $\Pi_j$, so bounding $\theta$ from below  amounts to bounding $\kappa$ away from $1$. Now put $\kappa=|\pi_i(y)|$ for a unitary vector $y\in \Pi_j$.  Then
$$
|\pi_i(y)|^2+|y-\pi_i(y)|^2=1,
$$
and since $y-\pi_i(y)\perp \ker u_i$, we get 
$$
|u_i||y-\pi_i(y)|\geq |u_i(y-\pi_i(y))|=|u_i(y)|.
$$
But since $u_j(y)=0$ and $|y|=1$ we have $|u_i(y)|\geq \eta$ by iii), so 
$$
\kappa^2=1-|y-\pi_i(y)|^2\leq 1-\frac{|u_i(y)|^2}{|u_i|^2}\leq 1-\frac{\eta^2}{|u_i|^2}.
$$ 
Finally the uniform bound $|u_i|\leq 2C$ yields the desired estimate (recall that $|s_i|_{\cc^1}\leq C$, while $|\eps_i|\ll1$). \cqfd
\subsection{Proof of theorem \ref{da}.}
In this paragraph, we review Donaldson and Auroux's works \cite{donaldson,donaldson2,auroux1,auroux2,auroux3} 
on the subject and indicate what must be changed to get theorem \ref{da}. 
Let us emphasize that our need for adapting these works mostly comes from the fact that the almost-complex structures
$J_j$ are not fixed. We must thus be very careful that the choices for $(\om_j,J_j)$ - which depend on $\eta$ as we saw 
above - do not affect the transversality estimates. This is not completely obvious because modifying $\om_j$ changes completely the line bundles in consideration, twisting them more and more when getting closer to $\om$. We claim 
however  that the decisive argument is already in Donaldson's original work : the estimates do not depend on the tensoring parameter $k$. The proof of theorem \ref{da} would however be much easier, would there exist an $\om$-tame almost 
complex structure $J$ whose set of compatible symplectic forms close to $\om$ projects to an open \nbd of $[\om]$ 
in $H^2(M,\R)$. Since we were unable to prove this point - which may well be false - we now proceed to a slightly fastidious adaptation.

We first explain the proof of theorem \ref{da} with only two bundles $(\cl_1,\cl_2)$ associated to $(\om_1,J_1)$, 
$(\om_2,J_2)$. That understood, the generalization to an arbitrary number $N$ of them will be straightforward. 
All Donaldson's construction relies on the existence of heavily localized approximately holomorphic sections. Namely, given 
$(M^{2n},\om,J)$ with $\om\in H^2(M,\Z)$ and $\cl$ a line-bundle on $M$ with connection of curvature $2i\pi\om$, 
Donaldson remarks :
\begin{lemma}\label{peak}
For all $p\in M$, there exists sections $\sigma_p^k$ of $\cl^{\otimes k}$ such that :
\begin{itemize}
\item[i)] $|\sigma_p^k(q)|\geq 1$ if $d_k(p,q)\leq 1$,
\item[ii)] $|\sigma_p^k(q)|_{\cc^1}\leq C_1e^{-C_2d_k(p,q)^2}$, 
\item[iii)] $|\bar \partial_J \sigma_p^k(q)|_{\cc^1}\leq \frac{C_1}{\sqrt k}e^{-C_2d_k(p,q)^2}$,
\item[iv)] the constants $C_1,C_2$ do not depend on $p$ nor $k$.
\end{itemize}
\end{lemma}
Usually the $k$ will be implicit and we denote $\sigma_p$ these sections. We must first check that this lemma can be extended to give sections $\sigma_{p,j}$ of $\cl_j^{\otimes k}$ ($j=1,2$) with the same estimates, where the constants $C_1,C_2$ are independant of the $(\om_j,J_j)$. This is possible because the dependancy of Darboux's theorem on the symplectic form can be made smooth. 
Indeed, the sections are of the form $\wdt \chi_k\circ f_k\circ \chi_k^{-1}(z)$ where 
\begin{itemize}
\item[\sbull] $\chi_k:\sqrt kB\hra M$ is the composition of the contraction of $\C^n$ $\delta_k:x\to x/\sqrt k$ and a Darboux chart 
$\chi_p:B(0,1)\hra M$  such that $\chi_p(0)=p$ and $\chi_p^*J(p)=i$,
\item[\sbull] $f_k:\sqrt kB \lra \C$ is (a far cut-off of) the map $f(z)=e^{-k|z|^2}$ viewed as a holomorphic section of the line bundle $\cl_\st^{\otimes k}$
with curvature $2i\pi k\om_\st$,
\item[\sbull] $\wdt\chi_k$ is a horizontal lift of $\chi_k$ to a bundle isomorphism between $(\cl_\st^{\otimes k},2i\pi k\om_\st)$
and $(\cl^{\otimes k},2i\pi k\om)$ above $\sqrt k B$.
\end{itemize}
All the estimates of lemma \ref{peak} come from the fact that $\chi_p$ can be chosen with a uniform bound on the derivatives. Since the uniform 
bound can be achieved not only with respect to $p$ but also with respect to the symplectic form in some \nbd of $\om$, the lemma 
holds for all $\om_j$ close enough to $\om$. According to \cite{donaldson}, the sections verifying estimates ii) and iii) above will be 
said approximately $J$-holomorphic.

The global construction of an  approximately holomorphic {\it and} uniformly transverse section is the following proposition (see \cite{donaldson, auroux3}) :
\begin{prop}\label{propagation} Given an approximately holomorphic sequence $(s_k)$ of sections of $\cl^{\otimes k}$, there exist points $(p_1,\dots,p_r)$ with $\cup B_k(p_i,1)=M$ and vectors  $(w^1,\dots,w^r)$ in $\C^{n+1}$ with $|w^i|\leq \delta$ such that 
the sequence
$$
s_w^k=s_k+\sum_{i=1}^r(w_0^i+\sum_{l=1}^n w^i_lz_l)\sigma_{p_i}^k,
$$   
is  approximately holomorphic and $\eta$-transverse, where $\eta$ does not depend on $k$ but only on $\delta$ and $g$ ($z_i$ denote the coordinates of the chart $\chi_{p_i}$).
\end{prop}
The number of  points involved in the process depends on $k$ and $g$ but on nothing else. This proposition 
relies itself on the following \cite{auroux3} :
\begin{thm}\label{yodoau}
Let $B^+:=B(\frac{11}{10})\subset \C^n$ and  $f:B^+\lra \C$. There exists $p\in \R$ depending only on $n$ such that if $|f|_{\cc^1(B^+)}\leq 1$ 
and $|\bar \partial f|_{\cc^1(B^+)}\leq \delta Q_p(\delta):=\delta |\ln \delta|^{-p}$, then there exists $w=(w_0,\dots,w_n)\in \C^{n+1}$ with 
$|w|<\delta$ and $f-w_0-\sum w_i z_i$ is $\delta Q_p(\delta)$-transverse to zero on $B(1)$.  
\end{thm}
The observation is now that since $\sigma_p(q)$ is large on $B^+(p_i):=B_k(p_i,\frac{11}{10})$ by  \ref{peak}, i), since 
$\bar \partial_{\wdt J}:=\bar\partial_{\chi_i^*J}$ is $k^{-\nicefrac{1}{2}}$-close to $\bar \partial$ on $\chi_i^{-1}(B^+(p_i))=B^+(0)\subset \C^n$ 
and since $\sigma_p$ and $(s_k)$ are approximately holomorphic,  we can apply the
 previous theorem to $f_i:=\nicefrac{s_k}{\sigma_{p_i}^k}$. This gives a $w$ for which $s_k-(w_0^i+\sum w_l^iz_l)\sigma_{p_i}^k$ is 
 $\delta Q_p(\delta)$-transverse to zero on $B_k(p_i,1)$. At this point again, provided that $J_j$ is close to $J$, the difference between $\bar\partial_{\wdt J_j}$ and $\bar \partial$ is of order $k^{-\nicefrac{1}{2}}$ with uniform constants and the argument extends to  our situation where $J_j$ is not fixed.
 
 The global construction then goes as follows. One can part the $r$ points into $K$ classes $\{\{{p_i}\}_{i\in I_\alpha}, \alpha\in [1,K]\}$ 
 and find constants $1>\delta_1>\dots >\delta_K$ with $\delta_{\alpha+1}=C\delta_\alpha Q_p(\delta_\alpha)$ such that : 
 \begin{itemize}
 \item[\sbull] The contributions of the $\{\sigma_{p_i}^k\}_{i\in I_\alpha}$ do not affect subsequently the transversality 
 at points of the same class. Precisely, points on a same class are sufficiently ($g_k$)-distant for the following to hold :  
 $$
 \begin{array}{l}
\left\{ \begin{array}{l}
 |w^i|\leq \delta_\alpha \\
 s_k \text{ is } 2\eta_\alpha=\delta_\alpha Q_p(\delta_\alpha)\text{-transverse on } B(p_{i'},1)
\end{array} \right\}
\Longrightarrow  \vspace*{,2cm}\\
\hspace*{1,5cm} \ds s_k+\sum_{\tiny\begin{array}{l}i\in I_\alpha\\ i\neq i' \end{array}}(w^i_0+\sum w_l^iz_l)\sigma_{p_i}^k 
\text{ is $\eta_\alpha$-transverse on } B(p_{i'},1).
\end{array}
 $$
 \item[\sbull] $C$ is a constant depending only on the constants $C_1,C_2$ of lemma \ref{peak} and $g$,  small enough that 
 the contributions of the $\{\sigma_{p_i}^k\}_{i\notin I_1\cup\dots\cup I_\alpha}$ does not affect the $\eta_\alpha$-transversality 
 on $V_\alpha:=\cup_{i\in I_1\cup\dots\cup I_\alpha} B_k(p_i,1)$. Precisely,
  $$
 \begin{array}{l}
\left\{ \begin{array}{l}
 |w^i|\leq \delta_{\alpha+1} \\
 s_k \; \eta_\alpha\text{-transverse on } V_\alpha
\end{array} \right\}
\Longrightarrow \vspace*{,2cm}\\
\hspace*{2cm}\ds s_k+\sum_{i\notin I_1\cup\dots\cup I_\alpha}(w_0^i+\sum w_l^iz_l)\sigma_{p_i}^k \text{ is 
$\frac{\eta_\alpha}{2}$-transverse on } V_\alpha.
\end{array}
 $$
 \item[\sbull] The number $K$ depends only on $C$, thus not on $k$ nor on $(\om_j,J_j)$.

 \end{itemize} 
Putting all this together, and using theorem \ref{yodoau} inductively on $B_k(p_i,1)$ for $i\in I_1,\dots,I_K$, we get proposition 
\ref{propagation} and thus Donaldson's theorem (starting with $(s_k)\equiv 0$). Since all constants $K,C,p$ above do not depend on $(\om_j,J_j)$ in a \nbd of $(\omega,J)$, we conclude that we can achieve the $\eta$-transversality with fixed 
$\eta$ ($=\eta_K/2$) for the sections $(s_j^k)$ of $\cl_j^{\otimes k}$ independantly of the approximation $\om_j$ we fixed.

We now give some details for the adaptation of the higher rank result, because it is the core of the difference (although nothing deep happens). The overall strategy is the same, but theorem \ref{yodoau} must be replaced by the following (see \cite{auroux3}) :
\begin{thm}\label{yodoau2}
Let $B^+:=B(\frac{11}{10})\in \C^n$ and  $f:B^+\lra \C^m$, $m\leq n$. There exists $p\in \R$ depending only on $n$ such that if $|f|_{\cc^1(B^+)}\leq 1$ 
and $|\bar \partial f|_{\cc^1(B^+)}\leq \delta Q_p(\delta)$, there exists $w=(w_0,\dots,w_n)\in \C^{m(n+1)}$ (each $w_i$ is a vector in $\C^m$) with  $|w|<\delta$ and $f-w_0-\sum w_i z_i$ is $\delta Q_p(\delta)$-transverse to zero on $B(1)$.  
\end{thm}
In order to apply it to our setting, decompose our section $s_k$ of $\cl_{1}^{\otimes k}\oplus \cl_{2}^{\otimes k}$ on $B_k(p,1)$ as 
$ s_k=(s^k_{1},s^k_{2})=f_1\sigma_{p,1}^k+f_2\sigma_{p,2}^k$ (identifying $\sigma_{p,1}^k$ with $(\sigma_{p,1}^k,0)$). 
The approximate holomorphicity of $s_k$ means that $|(\bar\partial_{J_1}s_1^k,\bar\partial_{J_2}s_2^k)|<Ck^{-\nicefrac{1}{2}}$, which implies in turn that 
$$
\left|\big(\bar\partial _{J_1}f_1,\bar\partial_{J_2}f_2\big)\right|<Ck^{-\frac{1}{2}} \hspace{1cm}\text{ on } B_k(p,1)
$$
because $\sigma_{p,1}^k$ and $\sigma_{p,2}^k$ are bounded below. In $\C^2$, we get 
$$
\left|\big(\bar\partial _{\wdt J_1}f_1\circ\chi_{p,1},\bar\partial_{\wdt J_2}f_2\circ\chi_{p,2}\big)\right|<Ck^{-\frac{1}{2}} \hspace{1cm}\text{ on } B_k(p,1).
$$
But $\wdt J_1,\wdt J_2$ are close to $i$ up to an order $k^{-\nicefrac{1}{2}}$ so $\bar\partial(f_1\circ\chi_{p,1},f_2\circ\chi_{p,2})$ is small. By theorem \ref{yodoau2}, we get a perturbation $(\wdt f_1,\wdt f_2)$ of $(f_1,f_2)$ given by theorem \ref{yodoau2} which is $\alpha$-transverse, 
{\it i.e.}
$$
\left|\big(\partial \wdt f_1\circ\chi_{p,1},\partial \wdt f_2\circ\chi_{p,2}\big)^{-1}\right|<\alpha \text{ whenever } |(\wdt f_1\circ\chi_{p,1},\wdt f_2\circ\chi_{p,2})|<\alpha.
$$ 
But again, since both $\partial_{\wdt J_1},\partial_{\wdt J_2}$ are $k^{-\nicefrac{1}{2}}$-close to the usual $\bar\partial$-operator, we get that for any point in $B_k(p,1)$ where $|(\wdt f_1,\wdt f_2)|<\alpha'$, $(\partial_{\wdt J_1}\wdt f_1\circ \chi_{p,1},\partial_{\wdt J_2}\wdt f_2\circ \chi_{p,2})=(\partial_{J_1}\wdt f_1,\partial_{J_2}\wdt f_2)$ has inverse of norm at most 
${\alpha'}^{-1}$, for $\alpha'$ slightly less that $\alpha$. Finally,
$$
(\partial_{J_1}\wdt s_1^k,\partial_{ J_2}\wdt s_2^k)=(\sigma_{1,p}\partial_{J_1}\wdt f_1,\sigma_{2,p}\partial_{J_2}\wdt f_2)+
(\wdt f_1\partial_{J_1}\sigma_{1,p},\wdt f_2\partial_{J_2}\sigma_{2,p}).
$$ 
Since $\sigma_j$ is bounded from below, $|(\sigma_{1,p}\partial_{J_1}\wdt f_1,\sigma_{2,p}\partial_{J_2}\wdt f_2)^{-1}|<(C\alpha')^{-1}$ 
(where $C$ is a universal constant), so for $|(\wdt f_1,\wdt f_2)|<\frac{C\alpha'}{2}$,
$$
\left| (\partial_{J_1}\wdt s_1^k,\partial_{ J_2}\wdt s_2^k)  \right|\leq \left(\frac{C\alpha'}{2}\right)^{-1}.
$$
This is the needed transversality for $(s_1,s_2)$. Getting it for all couples $(s_{j_1},s_{j_2})$ is then only a matter of induction over these couples, considering much smaller perturbations at each step. This is possible because we never destroy the 
approximate holomorphicity during this induction.
\cqfd    

\section{Desingularization and Biran decompositions.}\label{desingsec}
The aim of  this section is to use  theorem \ref{singdec} to give a generalization of Biran's decomposition's theorem to situations where the polarization is not smooth. Although nothing prevents a general study,  I prefer discussing an easy and concrete example in order to illustrate this point.

Consider $(\P^2,\om_\textnormal{FS})$ normalized so that the symplectic area of a projective line is $1$. Given our normalization of the standard form on $\R^{2n}$, this means that $\P^2$ is compactification of the ball of radius $1$. Any smooth cubic $C$ of $\P^2$ is a polarization of degree $3$, hence gives rise to an embedding of a standard disc bundle of radius $\nicefrac{1}{3}$ over $C$ by \cite{biran3} and to a full packing of $\P^2$ by one ellipsoid $\ce(3,\nicefrac{1}{3})$. The question 
studied in this paragraph is : what can we say when $C$ is a singular cubic of $\P^2$ instead of a smooth one ? As we shall see, although theorem \ref{singdec} does not formally consider singular curves, it can be easily associated to the classical desingularization techniques of algebraic geometry to provide a relevant answer to this question. 
\begin{thm} Let $C$ be a singular cubic of $\P^2$ with self-intersection at a point $p$. There exists a full packing of 
$\P^2$ by 
$$
B(\mu)\sqcup \ce(3-2\mu,\frac{1}{3})\sqcup \ce(\mu,\frac{2}{3}-\mu) \hspace{.5cm} \text{ for all } \mu<\frac{2}{3}.
$$
 Moreover, the cubic is covered by $B(\mu)$ - which it intersects along two Hopf discs, of area $\mu$ - and $\ce(3-2\mu,\nicefrac{1}{3})$ - which it intersects along the big axis, of area $3-2\mu$. It does not intersect $\ce(\mu,\nicefrac{2}{3}-\mu)$.
\end{thm}
 \noindent {\it Proof :} Assume for the moment that there exists a ball $B(\mu)$ centered at $p$ and whose intersection with $C$ is exactly two Hopf discs (this is certainly true for small $\mu$). Blowing-up this ball, we get the symplectic manifold $(\hat \P^2_1,\hat \om)$, where $[\hat \om]=l-\mu e$,
 endowed with a curve $\hat C$ (the strict transform of $C$) in the homology class of $3L-2E$.  The curves $\hat C$ and $E$ 
 are now smooth symplectic curves which intersect exactly twice, positively. They constitute a singular polarization of 
 $(\hat P^2_1,\hat \om)$ when $\mu<\nicefrac{2}{3}$, with 
 $$
 [\hat \om]=l-\mu e=\frac{1}{3}(3l-2e)+(\frac{2}{3}-\mu)e.
 $$
 By theorem \ref{singdec}, $(\hat \P^2_1;\hat \om)$ has a full packing by 
 $$
 \ce(3-2\mu,\nicefrac{1}{3})\sqcup \ce(\mu,\nicefrac{2}{3}-\mu).
 $$
 Moreover, it is easy to see that the disc $\{z_1=0\}$ can be brought out of $\ce(a,b)\subset \C^2$ by a symplectic isotopy with support in a small \nbd of $\ce(a,b)$. Thus, since $E\cap \ce(\mu-\eps,\nicefrac{2}{3}-\mu)$ is a Hopf disc and $E\cap \ce(3-2\mu-\eps,\nicefrac{1}{3})=\emptyset$, the manifold $\hat \P^2_1\priv E,\hat \om$ also has full packing by  $\ce(3-2\mu,\nicefrac{1}{3})\sqcup \ce(\mu,\nicefrac{2}{3}-\mu)$. 
 Blowing-down the exceptional divisor $E$ back, we therefore get a full packing of $\P^2$ as announced. The result for any $\mu$ is now a consequence of the next lemma.\cqfd

\begin{lemma} \label{cubichopf}
 For any singular cubic $C$ and for any $\mu<1$, there exists a ball $B(\mu)$ of capacity $\mu$ centered at the self-intersection point of $C$ and whose intersection with $C$ consists exactly of two Hopf discs.
 \end{lemma}
\noindent {\it Proof :} We show in fact that for any ball $B(\mu)$ there exists a cubic whose intersection with the ball is two Hopf discs. Since any two singular cubic are isotopic in $\P^2$, the lemma follows  \cite{barraud}. The proof is based on the blow-up construction of McDuff, which we do not review here (see for instance \cite{mcpo,mcsa}), and Gromov's theory of pseudo-holomorphic curves. Let $B(\mu)$ be a one-parameter family of balls of capacity $\mu$ in $\P^2$ for $\mu\in]0,1[$, $J_\mu$ an almost-complex structure suited for blow-up and $p_1(\mu),\dots,p_6(\mu)$ six generic points outside $B(\mu)$. Calling $p_0$ the center of $B(\mu)$, genericity means here that no three of the points $p_0,\dots,p_6$
lie in a same $J_\mu$-line and no six of them lie in a same $J_\mu$-conic. Denote also by $(\hat \P^2_1,\om_\mu,J_\mu)$ the symplectic blow-up of $\P^2$ along $B(\mu)$ endowed with the induced almost-complex structure (see \cite{mcsa}). Proving lemma \ref{cubichopf} then amounts to prove that the moduli space 
$$
\cm_\mu:=\{u:\P^1\lra\hat P^2_1\; |\; du\circ i=J_\mu\circ du , \; [u]=3L-2E,\; (p_1,\dots,p_6)\in \im u\},
$$
is not empty for $\mu<1$. We can also assume that the path of almost-complex structures $J_\mu$ is generic since we can modify $J_\mu$ in a \nbd of $p_1$ and $3L-2E$ is primitive. For $\mu$ small enough, this moduli space is obviously non-empty and even consists of exactly one point. If $\cm_\mu$ is empty for some $\mu$, there must be bubbling by Gromov's compactness theorem. This means that the class $3L-2E$ splits into a sum of classes $A_1+\dots+A_n$, where $A_i=k_iL-l_iE$ which are represented
by $J_\mu$-holomorphic curves. Since $E$ is also represented by a $J_\mu$-complex curve, we see by positivity of intersection that :
$$
(l_1,\dots,l_n)\in\{(1,1,0\dots,0),(2,0,\dots,0)\}.
$$
Moreover, $k_1>0$, as well as $k_2$ if $l_2=1$. Using now the positivity of intersection between the $A_i$,  the only possibilities are the following.
\begin{itemize}
\item[\sbull] All $k_i$ are positive (thus equal to one) if the decomposition consists of more that two terms $A_i$. Blowing down back to $\P^2$, the seven points must lie in a configuration of three lines, either two of them passing through  the center of $B(\mu)$, or one of them having a self-intersection point at $p_0$. But both configuration are impossible, one because a holomorphic line cannot have self-intersection, the other because the points were chosen generically. 
\item[\sbull] The decomposition is $(A_1,A_2)=(L-2E,2L)$ or $(2L-2E,L)$ if $k_1=2$. Blowing down again, we see that the first decomposition is impossible because it leads to a holomorphic line with a self-intersection point. The latter leads to a 
configuration of one line and one singular conic (meaning two lines passing through the center of the ball) passing through $(p_0,\dots,p_6)$, which is again impossible by genericity of the choice of the points.  
\item[\sbull] Or the decomposition is $(A_1,A_2)=(2L-E,L-E)$. Blowing down we get a configuration of one conic and one line 
intersecting at $p_0$, passing through a total of seven points, again impossible.  \cqfd
\end{itemize}

\section{Application to symplectic isotopies.}\label{isotopysec}
In \cite{moi5}, I explain a construction for isotopying balls. The principle is the following. Given a symplectic ball $B\subset (M^4,\om)$ (meaning that $B$ is the symplectic image of a 4-dimensional euclidean ball), define a supporting polarization for 
$B$ to be any {\it smooth } polarization $\Sigma$ of $M$ whose intersection with $B$ is exactly a Hopf disc in $B$ (the image of the intersection of $B^4\subset \C^2$ with a complex line). Very roughly, when there is a supporting polarization of degree $k$ for a ball of capacity less than $k^{-1}$, this ball can be brought into a standard position by symplectic isotopy. A precise statement is the following :
\begin{thm}\label{classicisotopy} Let $B_1,B_2\subset (M^4,\om)$ be symplectic balls of a rational symplectic manifold. Assume that :
\begin{itemize}
\item[\sbull] $B_1,B_2$ have supporting polarizations $\Sigma_1,\Sigma_2$ of same degree $k$,
\item[\sbull] $B_1,B_2$ have same capacity $c<k^{-1}$,
\item[\sbull] $\Sigma_1$ and $\Sigma_2$ are symplectic isotopic.
\end{itemize} 
Then $B_1,B_2$ are symplectic isotopic.
\end{thm}
The idea is that a given polarization allows to construct balls supported by this polarization in a very easy and flexible way. Conversely, any ball with this polarization as a supporting curve can be realized by such a construction. This theorem applies to some manifolds like $\P^2$ or $(S^2\times S^2,\om\oplus\om)$, but it is helpless for irrational symplectic manifolds, where there are no polarization at all. Even more unsatisfactory is the inaccuracy of the method for some very simple rational manifolds. For instance, when $\mu\in\Q\priv \Z$, the smooth polarizations of $(S^2\times S^2,\om\oplus \mu\om)$ have genus, and are therefore much more difficult to isotop, or even to bend to a supporting polarization, than spheres. As we will see below, 
this paper shows that singular polarizations are as good as smooth ones for the purpose of isotopies. This remark 
may be interesting in two respects. First, it sometimes allows to shortcutting any need for higher genus GW-invariants (for instance in the case of $S^2\times S^2$ as explained above). The second point is that singular polarizations may be in practice more {\it stable} 
objects than smooth ones, because they may arise as degeneracy of smooth polarizations through bubbling for instance. In view of the way the supporting polarization are produced (using pseudo-holomorphic curves), this stability property can be useful. We illustrate here the first point through an example :
\begin{thm}\label{singisotopy}
Any two balls of $(S^2\times S^2,\om\oplus\mu\om)$ are symplectic isotopic. 
\end{thm}
Below is a sketch of the proof. For more details, see also \cite{moi5} which is really devoted to the matter of isotopies. My aim here is only to explain how to use theorem \ref{singdec} in a problem of isotopy, when the method exposed in \cite{moi5} does not apply.

\noindent {\it Sketch of proof of theorem \ref{singisotopy} :} Let assume without loss of generality that $\mu>1$. Consider two symplectic balls $B_1,B_2$ of $M=(S^2\times S^2,\om\oplus\mu\om)$ of same capacity $c$ ($c<1$ by the non-squeezing theorem). By standard SFT arguments (stretching the neck) or blowing-up, it is easy to find supporting curves 
$\Sigma_i$ of $B_i$ in the homology class of $[S^2\times\{*\}]$ and symplectic curves $\Sigma_i'$ homological to $[\{*\}\times S^2]$ which do not meet $B_i$. Notice that $(\Sigma_i,\Sigma_i')$, $i=1,2$ are singular polarizations of $M$ in the sense of the present paper. Now by standard arguments, and because $\Sigma_i,\Sigma_i'$ are  spheres, the two couples of curves can be isotop one to another. The two balls can therefore be assumed to share a common singular supporting polarization. 
Notice now that a singular polarization $(\Sigma,\Sigma')$ gives rise to embeddings of an ellipsoid $\ce(1,\mu)$, which contains naturally a ball of capacity $c$, by paragraph \ref{embell}. As in the smooth case, these embeddings are completely determined by the single data of a Liouville form on $M\priv (\Sigma\cup\Sigma')$. The remaining of the reasonning is now exactly the same as in \cite{moi5} and we do not repeat it here : passing from $B_1$ to $B_2$ is only a matter of interpolating between two Liouville forms on $M\priv (\Sigma\cup\Sigma')$ which is easy.\cqfd

{\footnotesize
\bibliographystyle{abbrv}
\bibliography{bib3.bib}
}

\vspace{0cm}
\noindent Emmanuel Opshtein,\\
Institut de Recherche Mathématique Avancée \\
UMR 7501, Université de Strasbourg et CNRS \\
7 rue René Descartes \\
67000 Strasbourg, France\\
opshtein@math.unistra.fr
\end{document}